\documentclass[12pt]{amsart}
\usepackage{amsmath,amsthm,amsfonts,latexsym,amssymb}
\usepackage{graphicx}

\calclayout\evensidemargin .1in

\newtheorem{theorem}{Theorem}[section]
\newtheorem{lemma}[theorem]{Lemma}

\theoremstyle{definition}

\newtheorem{remark}[theorem]{Remark}

\title{A NOTE ON LEFSCHETZ FIBRATIONS ON COMPACT STEIN 4-MANIFOLDS}

\author{Selman Akbulut}
\address{Department of Mathematics, Michigan State University, Lansing MI, USA}
\email{akbulut@math.msu.edu}
\thanks{The first author is partially supported by NSF grant DMS 0905917}

\author{M. Firat Arikan}
\address{Department of Mathematics, University of Rochester, Rochester NY, USA}
\email{arikan@math.rochester.edu}
\date{\today}

\begin{document}

\begin{abstract}
Loi-Piergallini and Akbulut-Ozbagci showed that every compact Stein surface admits a Lefschetz fibration over the disk $D^2$ with bounded fibers. In this note we give a more intrinsic alternative proof of this result.
\end{abstract}

\maketitle


\section{Introduction}
In \cite{AO} (also \cite{LP} and \cite{P}) it was proven that every compact
Stein surface  admits a positive allowable Lefschetz
fibration over $D^2$ with bounded fibers (PALF in short), and conversely in \cite{AO} it was shown that every $4$-dimensional positive Lefschetz fibration over $D^2$ with bounded fibers is a Stein surface. The proof of \cite{AO} uses the fact that every torus link is fibered in $S^3$. Here we prove this by using a more intrinsic different approach, namely by an algorithmic use of positive stabilizations. This new approach is more closely related Giroux's proof of constructing open books to contact manifolds via ``contact cell decomposition" \cite{Gi}. The algorithm in \cite{A} constructs compatible open books for contact structures on $3$-manifolds using their surgery representations. The algorithm here is for $4$-manifolds, it has a similar technique and constructs PALF's on compact Stein surfaces starting from their handle diagrams which are explained briefly in \cite{Go}. We will give a different proof of the following result:

\begin{theorem} \label{thm:Main_Theorem}
Any compact Stein surface $W^4$ admits infinitely many pairwise inequivalent PALF's. Moreover, the corresponding open books on $\partial W$ supports the contact structure induced by the Stein structure on $W$.
\end{theorem}

We refer the reader to \cite{GS, K} for Lefschetz fibrations, to \cite{E, Go} for Stein manifolds, to \cite{Et1, Et2, Ge, Gi} for contact structures and
open books.

\section{An alternative proof of Theorem \ref{thm:Main_Theorem}}
\label{sec:Proof}

\begin{figure}[ht]
   \includegraphics{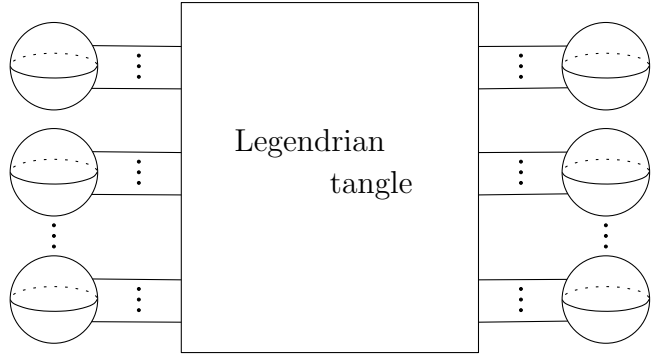}
   \caption{Handle diagram describing $W$}
  \label{fig:Figure_1}
\end{figure}

Let $W$ be a compact $4$-manifold admitting a Stein structure. By \cite{E}, $W$ has a handle decomposition which consists of a single $0$-handle,
1-handles, and $2$-handles attached to the union of the $0$-handle and the $1$-handles along some Legendrian knots $L_1, ..., L_n$ with framing $tb(L_i)-1$ where ``$tb$" denotes the Thurston-Bennequin framing. So by using \cite{Go} we can describe $W$ by a standard handle diagram given in Figure \ref{fig:Figure_1} where Legendrian tangle contains all crossings of the Legendrian link $L=L_1 \cup ...\cup L_n$ and for each crossing we assume the convention that the part of $L$ having more negative slope crosses over the one having less negative slope. \\

There are two cases that we need to consider:\\

\noindent \textbf{Case 1:} \emph{If there are no 1-handles in $W$.}\\

Suppose that $W$ is obtained from $D^4$ by attaching $2$-handles $H_1, ..., H_m$ along a Legendrian link $L=L_1 \cup ...\cup L_m$ sitting in $(S^3, \xi_{std})$. We will modify the algorithm of \cite{A} to construct a PALF on $D^4$ where we can realize each component of $L$ on a page of the PALF. Note that the algorithm of \cite{A} also guarantees that the page framing of $L_i$ is equal to $tb(L_i)$ for each $i$. Therefore, once we realize $L$ on pages of the PALF of $D^4$, attaching each $H_i$ along $L_i$ will extend the PALF structure, and we will be done.

\vspace{.05in}

\begin{figure}[ht]
   \includegraphics[width=.90\textwidth]{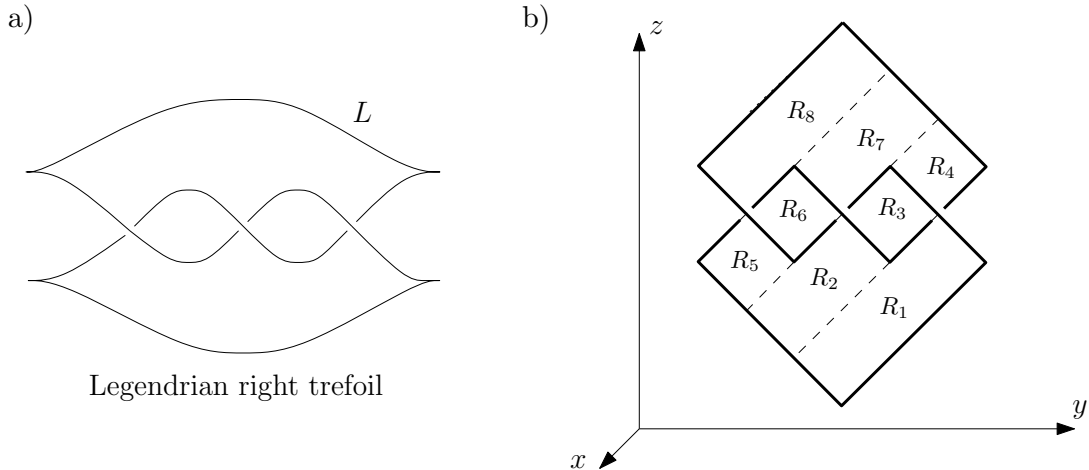}
   \caption{a) Legendrian $L \subset (S^3,\xi_{std})$, b) Front projection of $L$}
  \label{fig:Figure_2}
\end{figure}

For simplicity we'll take $L$ to be Legendrian right trefoil in our pictures. Given $L \subset (S^3,\xi_{std})$ we consider its front projection onto the $yz$-plane in $(\mathbb{R}^3, \xi_0=\textrm{Ker}(dz+xdy))$ as in Figure \ref{fig:Figure_2}. We divide the interior of the projection into rectangles $\{R_{i} \}$ using the lines with slope  $+1$ (see Figure \ref{fig:Figure_2}-b). Note that $i$ increases from down to up and right to left.

Here the main difference is that the Legendrian link $L$ (equipped with contact surgery coefficients) describes a contact surgery diagram (for some contact manifold) in \cite{A}, whereas here in our case it describes a compact Stein surface obtained from $D^4$.

\vspace{.05in}

\begin{figure}[ht]
  \begin{center}
   \includegraphics[width=.42\textwidth]{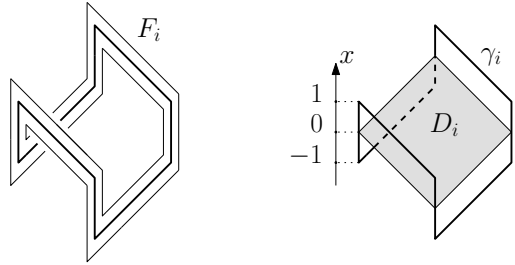}
   \caption{Constructing the Hopf band $F_i$ and the Legendrian unknot $\gamma_i$}
  \label{fig:Figure_3}
    \end{center}
\end{figure}

As in \cite{A}, for each $R_i$, we construct the Hopf band $F_i$ in $(\mathbb{R}^3, \xi_0) \subset (S^3,\xi_{std})$ by following the contact planes. Also we push the opposite sides of $R_i$ along the positive and negative $x$-axis and glue them using cords along the $x$-axis. This gives us a Legendrian unknot $\gamma_i$ sitting on $F_i$ with page framing equal to $tb(\gamma_i)=-1$ (see Figure \ref{fig:Figure_3}).

\begin{figure}[ht]
  \begin{center}
   \includegraphics[width=.65\textwidth]{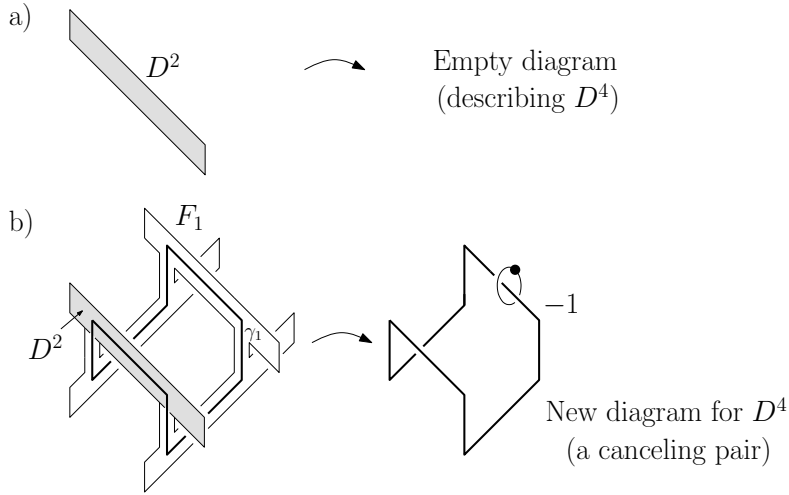}
   \caption{a) A regular fiber $D^2$ of the trivial PALF on $D^4$, b) a new PALF structure on $D^4$ with a regular fiber $F_1$}
  \label{fig:Figure_4}
    \end{center}
\end{figure}

We first consider the trivial PALF on $D^4$ with fibers $D^2$ and the trivial monodromy as in Figure \ref{fig:Figure_4}-a.
Note that the unique Stein structure on $D^4$ induces the unique tight structure $\xi_{std}$ on $S^3$, and this trivial PALF induces the compatible open book for $\xi_{std}$ with the same pages (fibers) and the monodromy. Now consider $R_1$ and glue the missing part of $F_1$ to the fiber $D^2$ and compose the positive Dehn twist $t_1$ along $\gamma_1$ with the existing (trivial) monodromy. This gives a new PALF structure on $D^4$ with a regular fiber $F_1$ and the monodromy $t_1$. In this process, we are actually positively stabilizing $D^4$ (see Figure \ref{fig:Figure_4}-b). Also the new corresponding open book on $S^3$ still supports $\xi_{std}$ by \cite{A}.

\newpage

Next we introduce other rectangles one by one (in ascending order) to the projection. For each $R_i$ introduced, we positively stabilize $D^4$ and extend the PALF structure as explained in the following remark.

\begin{remark} \label{rem:extending_PALF}
To be able to extend the PALF structure, we must introduce the rectangles in a special order. Such order is provided by how we number the rectangles above. More precisely, the above ordering guaranties that attaching the missing part of the Hopf band $F_i$ to the fiber, say $S_{i-1}$, of the PALF corresponding to $R_1,...R_{i-1}$ is equivalent to plumbing a positive Hopf band to $S_{i-1}$, and so the new surface $S_{i}$ is a fiber of a new PALF. We will explicitly show this equivalence only for Case 2 below (see Lemma \ref{lem:extending_PALFs}). Showing the equivalence for Case 1 is straightforward (comparing to that for Case 2). Therefore, we leave it to the reader.
\end{remark}

\begin{figure}[ht]
   \includegraphics[width=.55\textwidth]{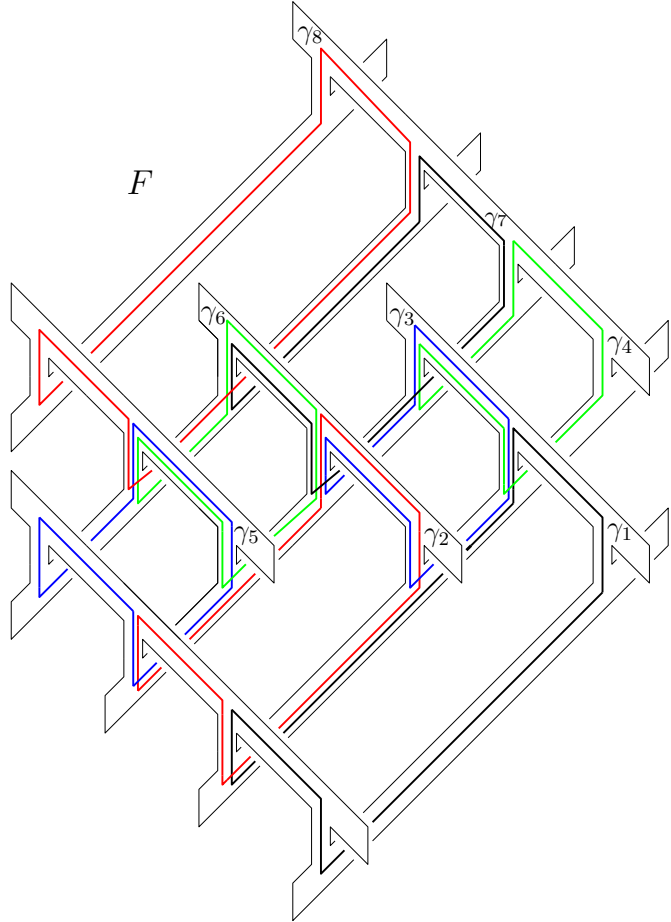}
   \caption{The fiber $F$ of the PALF on $D^4$ containing $L$}
  \label{fig:Figure_5}
\end{figure}

\begin{figure}[ht]
   \includegraphics[width=.45\textwidth]{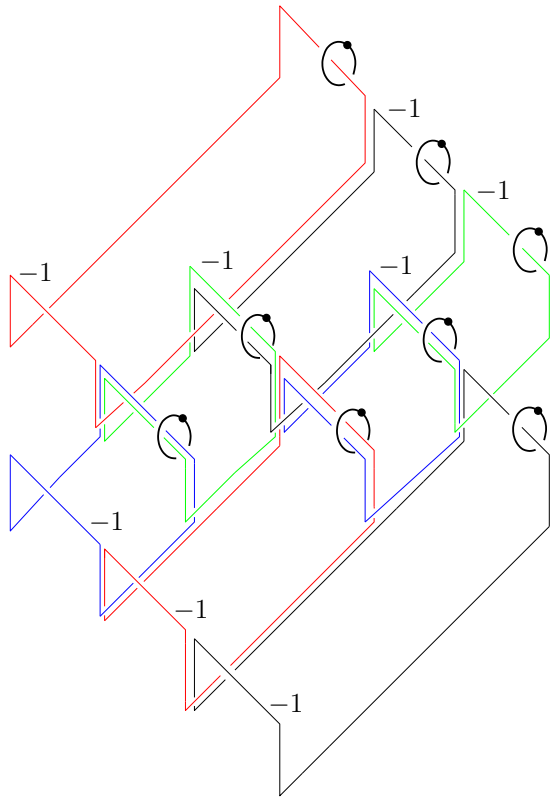}
   \caption{The diagram for $D^4$ corresponding to the PALF in Figure \ref{fig:Figure_5} $\quad$ ($n$ canceling pairs; they are canceled from right to left)}
  \label{fig:Figure_6}
\end{figure}

When all rectangles are introduced, we get a new PALF structure on $D^4$ with bounded fiber $F$ (union of Hopf bands $F_1,...,F_n$) and monodromy $t_1t_2...t_n$ where each $t_i$ is the positive Dehn twist along the Legendrian unknot $\gamma_i$ coming from the stabilization corresponding to $R_i$ (see Figure \ref{fig:Figure_5}). Note that this process does not change $D^4$ because in the corresponding handle diagram we have $n$ canceling pairs of $1$- and $2$-handles as shown in Figure \ref{fig:Figure_6} where we use the convention that $\gamma_i$ crosses over $\gamma_j$ if $i>j$. Here we consider each $\gamma$-curve in a different page $F$ of the corresponding open book supporting $\xi_{std}$ by pushing them in the (pointing out) normal direction of $F$.

\begin{remark} \label{rem:isotopic_strs_links}
We note that the induced open book on $S^3$ has the monodromy $t_1t_2...t_n$ (product of positive Dehn twists), and so, by the uniqueness, it is compatible with a contact structure isotopic to $\xi_{std}$. Also we have obtained a Legendrian link on a page $F$ which is topologically equivalent to $L$. By following similar arguments in the appendix section of \cite{P}, we immediately see that this new Legendrian link is, indeed, Legendrian isotopic to the original link $L$. For simplicity, we denote the new link also by $L$.
\end{remark}

So far we have constructed a PALF structure on $D^4$ such that the Legendrian link $L=L_1 \cup ... \cup L_m$ is embedded on a page $F$ of the open book supporting $\xi_{std}$. Also the framing of each $L_i$ coming from $F$ is equal to $tb(L_i)$.
Therefore, when we attach $2$-handles $H_1,...,H_m$ along the Legendrian knots $L_1,...,L_m$, we do not only get the compact Stein surface $W$ but also extend the PALF on $D^4$ to a PALF on $W$. Note that a regular fiber of the resulting PALF on $W$ is still $F$, but the monodromy is now equal to $$t_1t_2...t_ns_1s_2...s_m$$ where $s_i$ is the positive Dehn twist along $L_i$ for each $i=1,...,m$.\\

\newpage

\noindent \textbf{Case 2:} \emph{If there are 1-handles in $W$.}\\

\begin{figure}[ht]
   \includegraphics{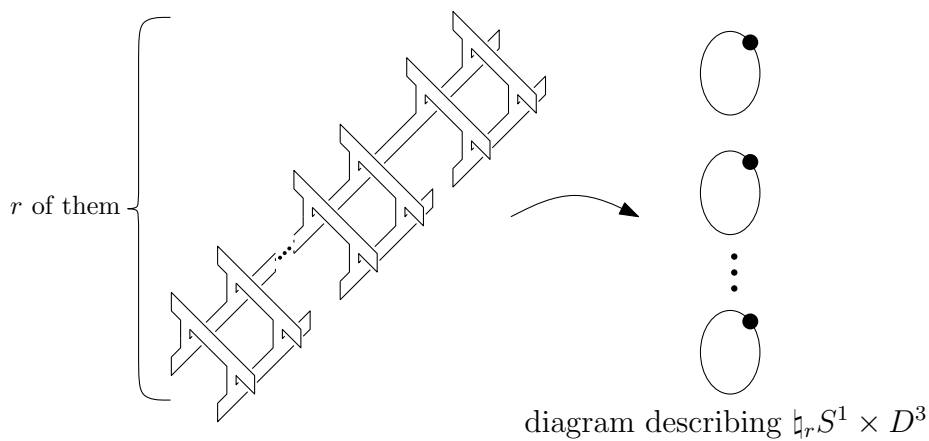}
   \caption{A fiber of the trivial PALF on $\natural_r S^1 \times D^3$ and the corresponding diagram}
  \label{fig:Figure_7}
\end{figure}

Suppose that $W$ is obtained from $D^4$ by attaching $r$ $1$-handles and $2$-handles $H_1, ..., H_m$ along a Legendrian link $L=L_1 \cup ...\cup L_m$ sitting in $(\#_r S^1 \times S^2, \eta_{std})$. A standard handle diagram for $W$ is given in Figure \ref{fig:Figure_1}. The union of $D^4$ and $1$-handles gives $\natural_r S^1 \times D^3$ whose Stein structure induces the unique tight structure on $\eta_{std}$ on $\#_r S^1 \times S^2$. First consider the trivial PALF (with trivial monodromy) on $\natural_r S^1 \times D^3$. A regular fiber of this PALF is given in Figure \ref{fig:Figure_7} where the reader should realize that the core circle of each Hopf band links the corresponding dotted circle once. Note that the corresponding open book supports $\eta_{std}$ and in the standard handle diagram of $W$, if a knot passing over a particular $1$-handle, then it must link to the core circle of the corresponding Hopf band once.

\begin{figure}[ht]
   \includegraphics{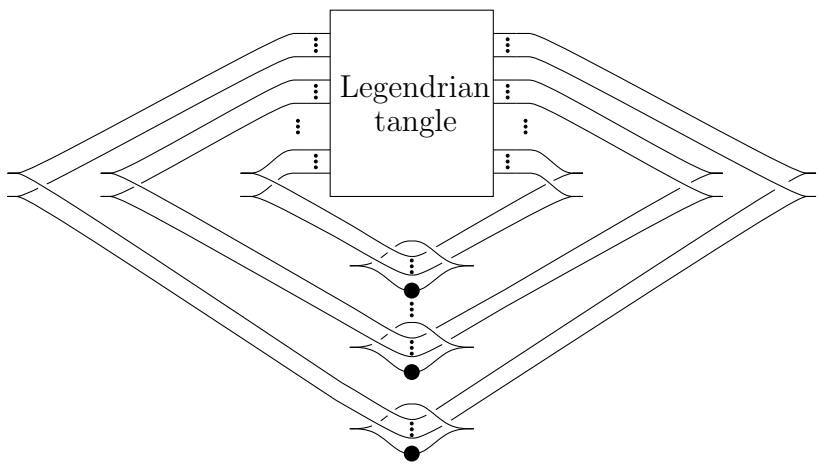}
   \caption{Converting the diagram}
  \label{fig:Figure_8}
\end{figure}

We first modify the handle diagram in Figure \ref{fig:Figure_1} by twisting the strands
going through each $1$-handle and replace $1$-handles with dotted Legendrian unknots as illustrated in Figure \ref{fig:Figure_8} (compare with \cite{AO}).

\begin{figure}[ht]
   \includegraphics{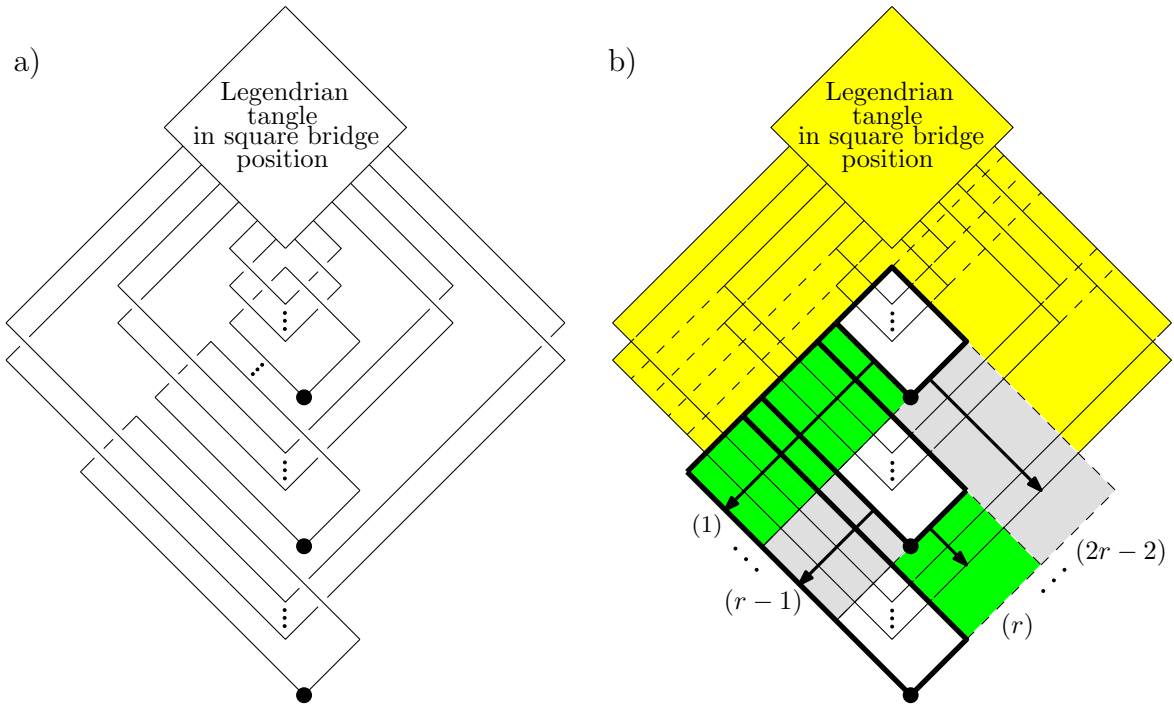}
   \caption{a) Projecting the diagram, b) Introducing and ordering the rectangles}
  \label{fig:Figure_9}
\end{figure}

Then by pretending the resulting diagram sits in $(\mathbb{R}^3,\xi_0)$, we consider its projection onto the $yz$-plane as in Figure \ref{fig:Figure_9}-a. Next we want to divide the interior of the projection into rectangles. Surely this can be done in many different ways. However, with a little  care we can decrease the number $n$ of rectangles as illustrated in Figure \ref{fig:Figure_9}-b where we assume that the bold arcs are introduced first and they correspond the fiber of the trivial PALF on $\natural_r S^1 \times D^3$ given in Figure \ref{fig:Figure_7}. Note that the region bounded by the bold arcs are divided into concentric rectangles, and that we add some additional arcs to the projection to be able to extend PALF structures (why we need these additional arcs will be explained below). We also number the rectangles as follows: The concentric rectangles inside bold squares come first. The squares in the colored regions $(1),(2),..., (2r-2)$ come second in the ascending order. (Here, for each of these regions, we number the rectangles in the order indicated by the arrow. Also, for the regions $(r),..., (2r-2)$, the boundary rectangles come first.) Finally, the rectangles in the yellow region come last in the order explained in Case 1.

As pointed out in Remark \ref{rem:extending_PALF}, we now explain how the ordering we choose in Figure \ref{fig:Figure_9}-b extends PALF structures.

Let $P$ be a diagram in the $yz$-plane divided into rectangles whose sides are on the lines of slopes $\pm 1$. Let $F_P$ denote the surface obtained by following the contact planes in $(\mathbb{R}^3, \xi_0)$ with the front projection $P$ (here we just generalize the construction of the positive Hopf band from a rectangle, see Figure \ref{fig:Figure_3}).

\begin{lemma} \label{lem:extending_PALFs}
In each case below, suppose that the surface $F_P$ is a fiber of a PALF. Then so is $F_{Q}$ where $Q$ is obtained from $P$ by adding the blue arc.
\end{lemma}

  \begin{center}
   \includegraphics{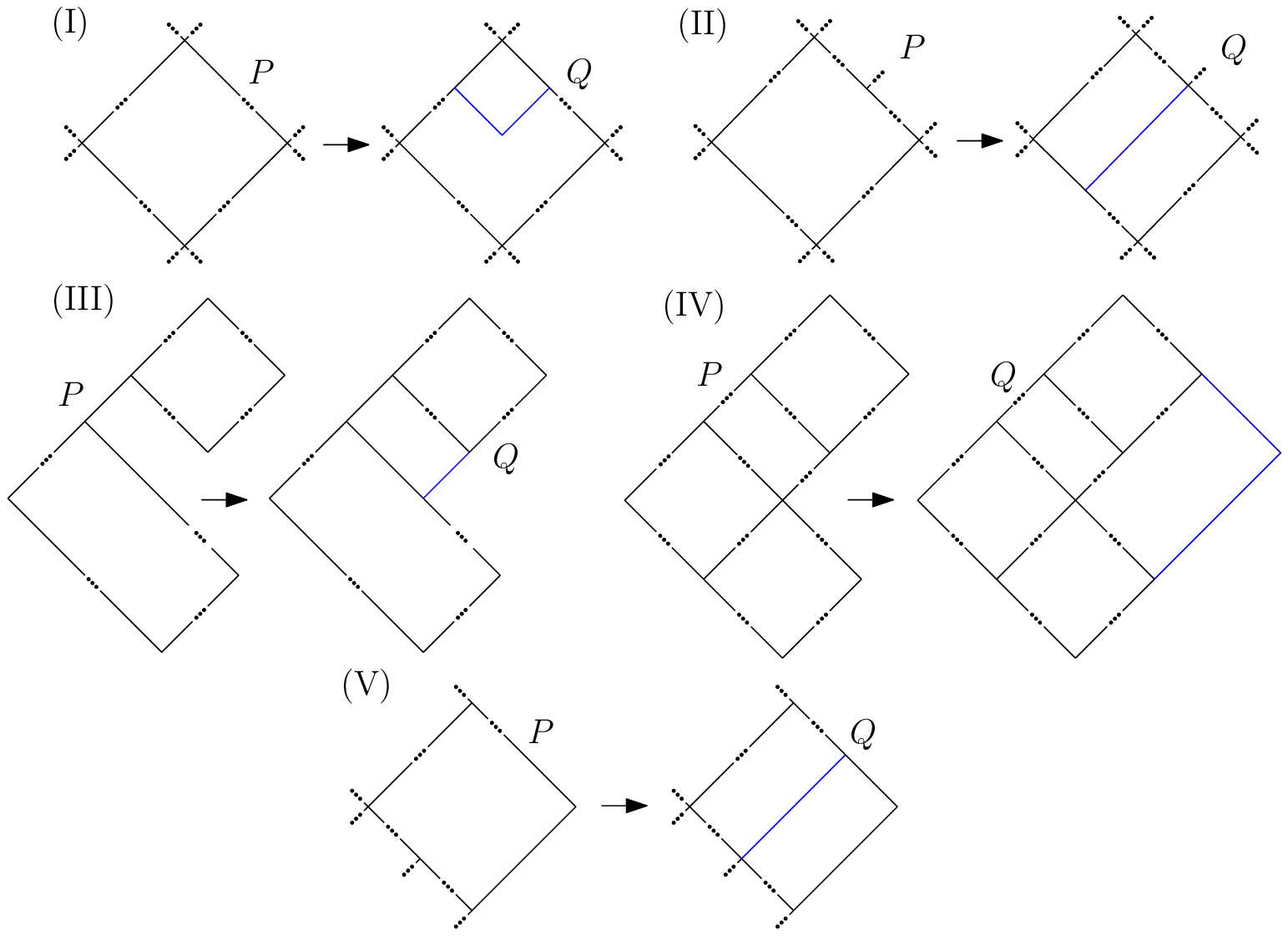}
     \end{center}

\begin{proof} In each case, the surface $F_Q$ is obtained from $F_P$ by adding the strip which projects onto the blue arc in the diagram $Q$. For each case, we draw a picture below where this strip is given in blue. In each picture, we show that adding the blue strip to the fiber $F_P$ is equivalent to plumbing a positive Hopf band $H^+$ to $F_P$ along the red arc in $F_P$. (Dashed arrows indicate the isotopies taking $H^+$ to the blue strip.) Thus, the surface $F_Q$ is a fiber of a PALF.
\end{proof}
\begin{figure}[ht]
   \includegraphics{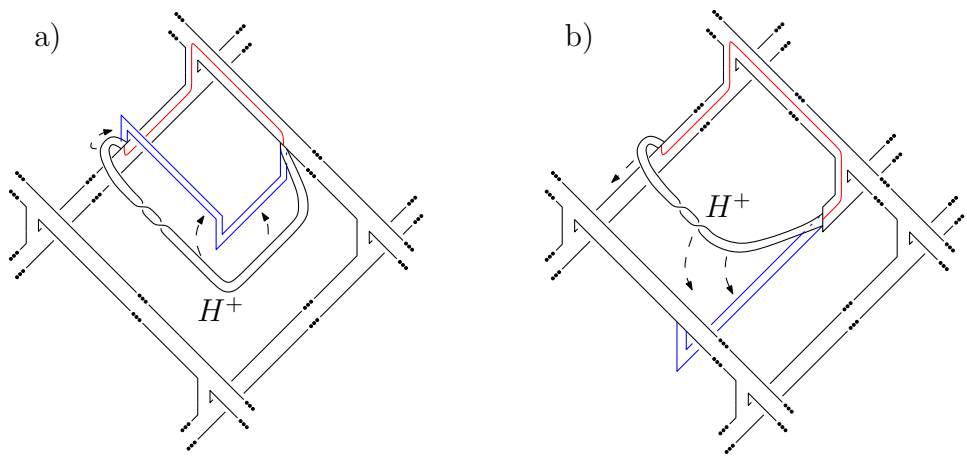}
   \caption{Proof of a) Lemma \ref{lem:extending_PALFs}-(I), b) Lemma \ref{lem:extending_PALFs}-(II)}
  \label{fig:Figure_10}
\end{figure}

\clearpage

\begin{figure}[ht]
   \includegraphics{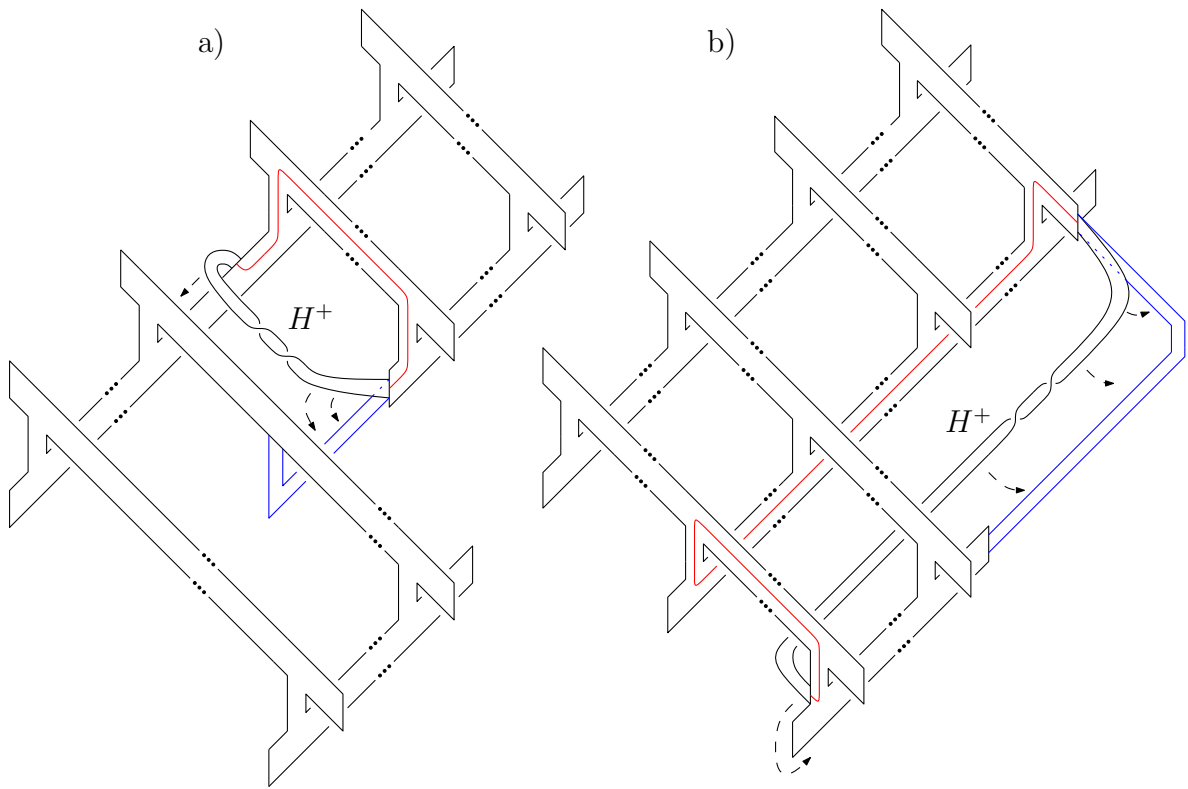}
   \caption{Proof of a) Lemma \ref{lem:extending_PALFs}-(III), b) Lemma \ref{lem:extending_PALFs}-(IV)}
  \label{fig:Figure_11}
\end{figure}

\begin{figure}[ht]
   \includegraphics{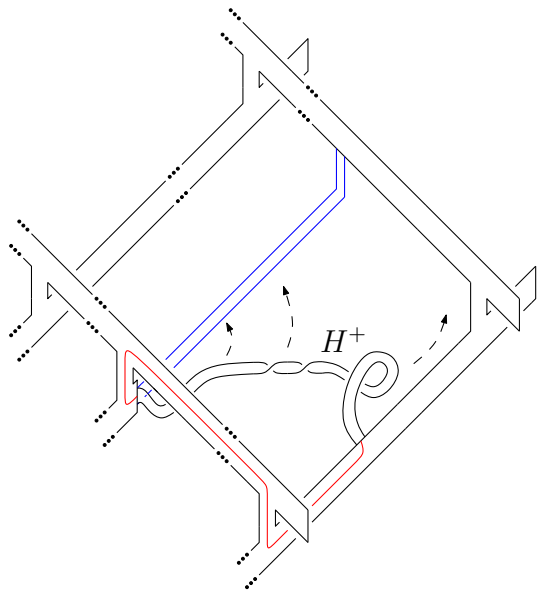}
   \caption{Proof of Lemma \ref{lem:extending_PALFs}-(V)}
  \label{fig:Figure_12}
\end{figure}

\clearpage

Now we introduce the rectangles in the described order to the diagram in Figure \ref{fig:Figure_9}-b, and accordingly positively stabilize $\natural_r S^1 \times D^3$. By Lemma \ref{lem:extending_PALFs}, we know that PALF structure extends for the rectangles in the regions $(1),..., (2r-2)$ and the bold rectangles. For the remaining rectangles (the ones in the yellow region in Figure \ref{fig:Figure_9}-b), the PALF structures extend as in Case 1. When all rectangles are introduced, we get a new PALF structure on $\natural_r S^1 \times D^3$ with bounded fiber $F$ and monodromy $T$ which is a product of (+)-Dehn twists $t_1,t_2,..., t_n$ in a certain order determined by the plumbings (here each $t_i$ is as before).

\begin{remark} As in Remark \ref{rem:isotopic_strs_links}, the monodromy of the induced open book on $\#_r S^1 \times S^2$ is the product of positive Dehn twists which implies, by the uniqueness, that it is compatible with a contact structure, say $\eta$, isotopic to $\eta_{std}$. Moreover, we have obtained a Legendrian link, topologically equivalent to $L=L_1 \cup ... \cup L_m$, on a page $F$ of the open book supporting $\eta_{std}$. This new Legendrian link is, indeed, Legendrian isotopic to the original link $L$ and its Legendrian type is the same with respect to $\eta$ and $\eta_{std}$. To see this, we modify the arguments in the appendix section of \cite{P} as follows: First, the binding $\partial F$ is not a $(p,q)$-torus knot but it is still a fibered link (in $\#_r S^1 \times S^2$), and so we can still consider the handlebodies $H_1, H_2$ (by thickening $F$ and taking its complement) and then apply Thurston-Winkelnkemper construction (as in the proof of Proposition 4 of \cite{P}) to get a new contact form which defines the contact structure $\eta$ coinciding $\eta_{std}$ along $H_1$. This shows that the new link on $F$ (which will replace $L$) is still Legendrian with respect to $\eta$. To check that its Legendrian type is the same with respect to $\eta$ and $\eta_{std}$ we show that in the handlebody $H_2$ the contact structures  $\eta$ and $\eta_{std}$ are isotopic relative to the boundary by modifying Lemma 3 in \cite{P} to our case: Note that for the subsurface $F'$ of $F$ constructed by plumping positive Hopf bands, the corresponding handlebody $H' \subset H_2$ can be thought of as solid tori with convex boundary and dividing set given by two parallel curves with slope $-1$, and the required isotopy between $\eta$ and $\eta_{std}$ follows by the induction argument in \cite{P}. However, we need to see that they are isotopic in the whole $H_2$. To this end, we decompose $H_2$ so that it is the union of $H'$ and a collection of $3$-balls. A typical $3$-ball is of the form $B=D^2 \times[0,1]$. We may assume that each cutting disk $D_i=D^2 \times \{i\} \subset H_2$ ($i=0,1$) has a Legendrian boundary $\partial D_i \subset \partial H_2$, and each $\partial D_i$ intersects the dividing set of $H'$ in two distinct points, say $p_i, q_i$. Now, since $B$ lives in a manifold with tight structure, the dividing set of the boundary $2$-sphere $\partial B$ consists of a single closed curve $C$ (otherwise contact structure would be overtwisted). By a small perturbation, we can assume that $S$ is convex. Also as $S$ and $\partial H'$ have the cutting disks $D_1,D_2$ in common, $C \cap D_i$ connects $p_i$ and $q_i$. Now using Lemma 3.11 (edge-rounding) in \cite{H}), we can glue $\overline{S-(D_1 \cup D_2)}$ and $\overline{\partial H -(D_1 \cup D_2)}$ to get $\partial H_2$ back. But this time we are sure that $\partial H_2$ is a convex surface with dividing set given by two parallel curves with slope $-1$. Finally, the required isotopy follows as explained in \cite{P}.
\end{remark}

As before we denote the new link on $F$ also by $L$. By the construction the framing of each $L_i$ coming from $F$ is equal to $tb(L_i)$. Therefore, as in Case 1 attaching $2$-handles $H_1,...,H_m$ along $L_1,...,L_m$ extends the PALF on $\natural_r S^1 \times D^3$ to a PALF on $W$. The resulting PALF on $W$ has a bounded regular fiber $F$, and its monodromy is $$T \cdot s_1s_2...s_m$$ where $s_i$ is the positive Dehn twist along $L_i$ for each $i=1,...,m$. \\

We note that in both of the above cases the final open book corresponding to the final PALF supports the contact structure on $\partial W$ induced by the Stein structure on $W$. The reader is referred to \cite{Et1} for details on compatibility. Also observe that once we have a PALF structure on $W$, we can get infinitely many pairwise inequivalent PALF's on $W$ by positively stabilizing the original one.\\

As a final remark, one should mention that all the steps in our algorithmic proof are applicable to any Stein manifold described by any standard handle diagram. To see more rigorous explanations the reader is referred to \cite{A}. $\square$

\section{Example}
\label{sec:Example}

\begin{figure}[ht]
  \begin{center}
   \includegraphics{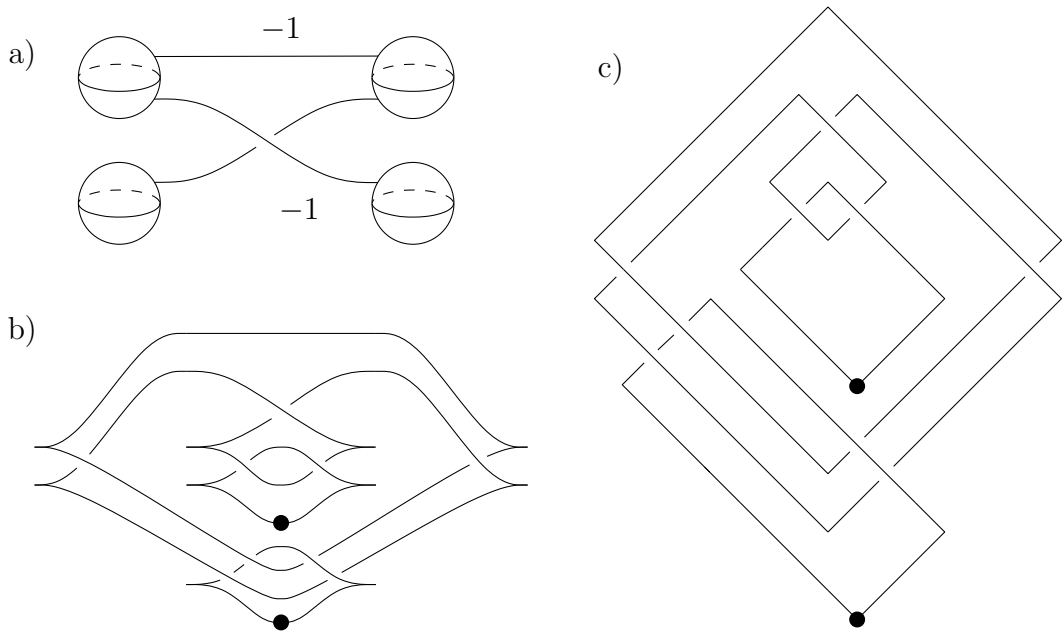}
   \caption{a) A Stein surface $W$ ($r=m=2$, coefficients are relative to $tb$-framing), b) converting the diagram, c) projecting onto the $yz$-plane}
  \label{fig:Figure_13}
    \end{center}
\end{figure}

As an example, we will apply our algorithm to construct a PALF structure on the compact Stein surface $W$ given in Figure \ref{fig:Figure_13}-a.

\begin{figure}[ht]
  \begin{center}
   \includegraphics{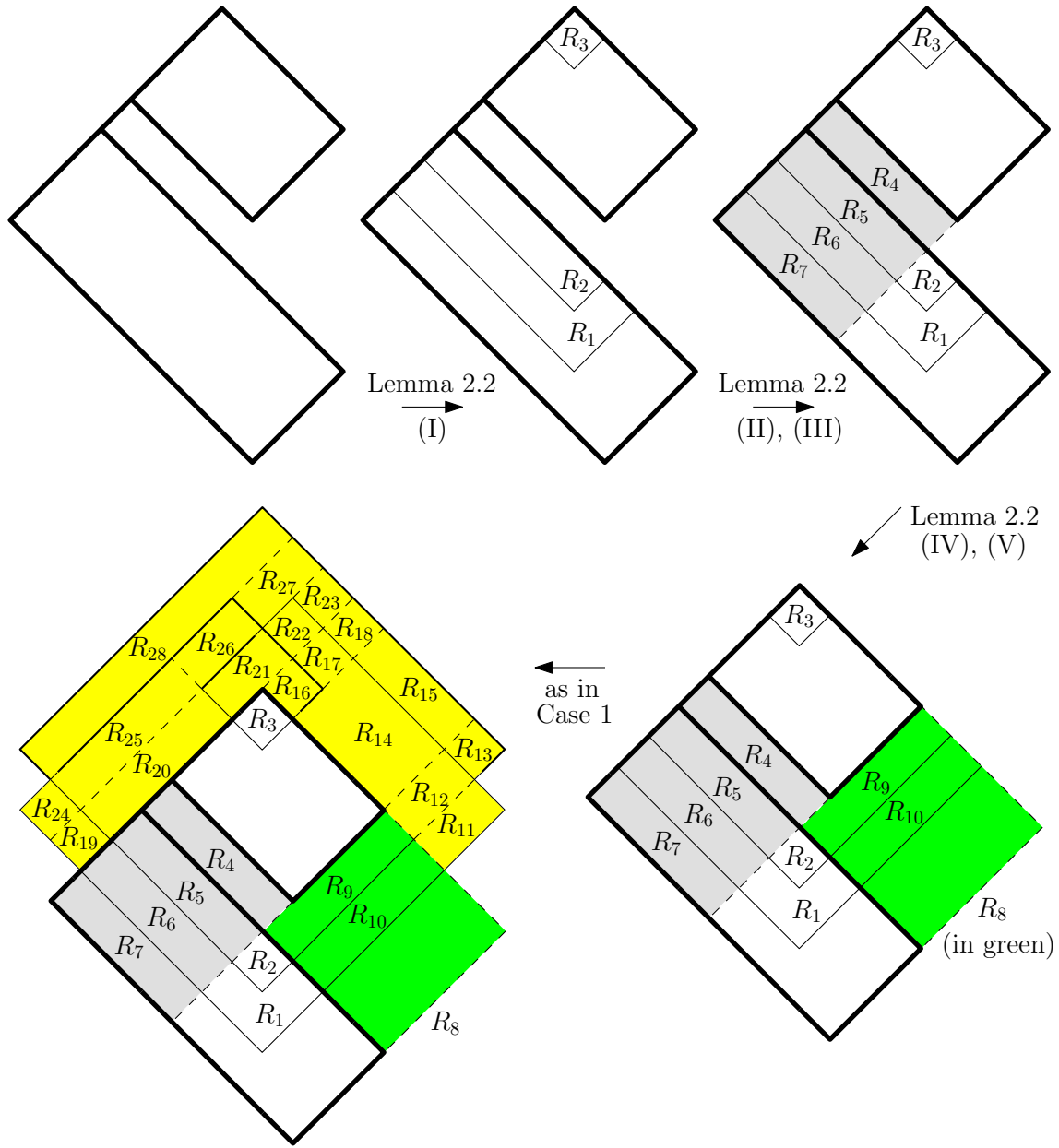}
   \caption{Introducing the rectangles in the order given by the algorithm.}
  \label{fig:Figure_14}
    \end{center}
\end{figure}

We first convert the $1$-handles into dotted circles and obtain the diagram in Figure \ref{fig:Figure_13}-b. Then we consider its projection onto the $yz$-plane as in Figure \ref{fig:Figure_13}-c. Next we introduce the rectangles to the projection in the order depicted in Figure \ref{fig:Figure_14}, and accordingly positively stabilize $\natural_r S^1 \times D^3$. By Lemma \ref{lem:extending_PALFs}, we know that PALF structure extends for the rectangles $R_1,..., R_{10}$. For $R_{11},...R_{28}$, we extend the PALF structures as in Case 1.\\

\begin{figure}[ht]
   \includegraphics{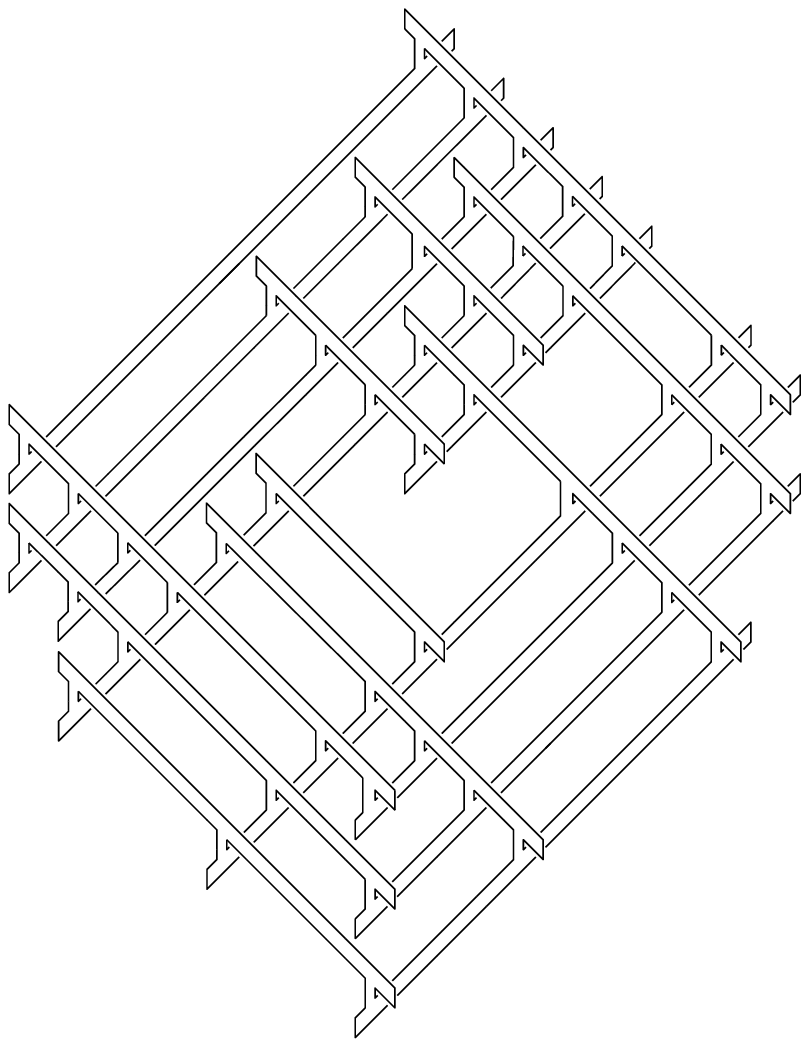}
   \caption{The fiber $F$ of the PALF on $\natural_2 S^1 \times D^3$}
  \label{fig:Figure_15}
\end{figure}

The fiber of the resulting PALF on $\natural_2 S^1 \times D^3$ is shown in Figure \ref{fig:Figure_15}. We remark that in the corresponding handle diagram we have $28$ canceling pairs of $1$- and $2$-handles.\\

Finally, note that the Legendrian link describing $W$ is now embedded on a page of the PALF on $\natural_2 S^1 \times D^3$ with page framing equal to tb-framing. Therefore, by attaching the $2$-handles, we obtain a PALF structure on $W$.

\clearpage

\addcontentsline{toc}{chapter}{\textsc{References}}

\end{document}